\documentclass[11pt]{article}

\usepackage{amsmath,amsthm,latexsym}
\usepackage[frame]{xy}

\newtheorem{theorem}{Theorem}[section]
\newtheorem{cor}[theorem]{Corollary}
\newtheorem{lemma}[theorem]{Lemma}
\newtheorem{remark}[theorem]{Remark}

\newtheorem{definition}[theorem]{Definition}

\newtheorem{prop}[theorem]{Proposition}
\makeatletter
\newfont{\footsc}{cmcsc10 at 8truept}
\newfont{\footbf}{cmbx10 at 8truept}
\newfont{\footrm}{cmr10 at 10truept}
\makeatother
\def\graphminus{-}

\def\schn{scheme}

\def\Dsch{Duchet scheme }
\def\Dschn{Duchet scheme}
\def\Dschs{Duchet schemes }
\def\hDsch{$H$-Duchet scheme }
\def\hDschn{$H$-Duchet scheme}

\def\Ssch{scheme }
\def\Sschn{scheme}
\def\Sschs{schemes }
\def\Shsch{$H$-scheme }
\def\Shschn{$H$-scheme}
\def\Shschs{$H$-schemes }
\def\Shpsch{$H'$-scheme }
\def\Shpschn{$H'$-scheme}

\def\Csch{colored scheme }
\def\Cschn{colored scheme}
\def\Cschs{colored schemes }
\def\Chsch{colored $H$-scheme }
\def\Chschn{colored $H$-scheme}

\def\rc{contractible }
\def\rcn{contractible}

\def\mpcontn{$M'$-contraction}
\def\mpc{$M'$-contractible }
\def\mpcn{$M'$-contractible}

\begin{document}

\title{Finding minors in graphs with a given path structure}
\author{
Andr\'e K\"undgen \\
\small  Department of Mathematics\\
\small Calfornia State University San Marcos\\
\small San Marcos, CA 92096\\
{\em  akundgen{@}csusm.edu}
\and
Michael J. Pelsmajer\thanks{The second author gratefully acknowledges the support from NSA Grant H98230-08-1-0043 and the Swiss National Science Foundation Grant No. 200021-125287/1.}  \\
\small Department of Applied Mathematics\\
\small Illinois Institute of Technology\\
\small Chicago, IL 60616\\
{\em  pelsmajer{@}iit.edu}
\and
Radhika Ramamurthi \\
\small  Department of Mathematics\\
\small Calfornia State University San Marcos \\
\small San Marcos, CA 92096\\
}

\maketitle
\bigskip\bigskip
\begin{abstract}
 Given graphs $G$ and $H$ with $V(H)\subseteq V(G)$, suppose that we have a $u,v$-path $P_{uv}$ in $G$ for each edge $uv$ in $H$.  There are obvious
 additional conditions that ensure that $G$ contains $H$ as a rooted
 subgraph, subdivision, or immersion; we seek conditions that
 ensure that $G$ contains $H$ as a rooted minor or minor.
 This naturally leads to studying sets of paths that form an $H$-immersion,
 with the additional property that paths that contain the same vertex must have a common
 endpoint.  We say that $H$ is {\em \rcn} if,
 whenever $G$ contains such an $H$-immersion, $G$ must also contain a rooted $H$-minor.
 We show, for example, that forests, cycles, $K_4$, and $K_{1,1,3}$ are \rcn, but that
 graphs that are not 6-colorable and graphs that contain certain subdivisions of $K_{2,3}$ are not \rcn.
\end{abstract}


\section{Introduction}

A common question in graph theory is:
Given two graphs $G$ and $H$, with $V(H) \subseteq V(G)$, how can $H$ be
represented within $G$?  The ``model'' of $H$ that we seek in $G$ may be a
(rooted) subgraph, subdivision, minor, or some other substructure.  Via
connectivity or structural arguments, one may be able to find a $u,v$-path $P_{uv}$ in
$G$ that ``models'' each edge $uv\in E(H)$; then the question is what
additional conditions will ensure that these paths yield the desired (rooted)
$H$-substructure.

Since additional structure only makes it easier to find $H$ represented within $G$,
we usually assume that $G$ is the {\em underlying graph} of the paths, 
i.e., $V(G)=\bigcup_{uv\in E(G)} V(P_{uv})$ and $E(G)=\bigcup_{uv\in E(G)} E(P_{uv})$.
Then, to guarantee a (rooted) $H$-subgraph, $H$-subdivision, and
$H$-immersion, there are obvious necessary and sufficient conditions; namely,
that the paths are length~$1$, internally vertex-disjoint, and internally
edge-disjoint, respectively.
The situation is less obvious for minors.

One obvious necessary condition for having $H$ as a rooted minor
is that the paths should not have any internal vertices in $V(H)$.
Another issue arises when $G$ contains a star or subdivision of a star $S$
with leaf set $V(H)$.  Then every edge $uv\in E(H)$ can be represented by a
path in $S$, but $G$ need not contain an $H$-minor, and the structure of $S$
contains no structural information about $H$.  To avoid this sort of thing,
we limit the variety of paths that can pass through the same vertex of $G$:

\begin{definition}\label{d:scheme-def}
Let $G$ and $H$ be graphs with $V(H)\subseteq V(G)$. A {\em \Shschn} is a collection of paths
$\{P_{uv}: uv\in E(H)\}$ (called {\em \schn-paths}) such that
\begin{enumerate}
\item $P_{uv}$ is a $u,v$-path in $G$ with  $V(P_{uv})\cap V(H)=\{u,v\}$, and
\item every collection of edges $e_1,\dots, e_k\in E(H)$ such that $V(P_{e_1})\cap \dots \cap V(P_{e_k})\neq \emptyset$ must have a common endpoint $u\in V(H)$,
i.e. $u$ is incident with each of $e_1,\dots e_k$.
\end{enumerate}
\end{definition}

In this paper, we are primarily interested in whether the presence of an \Shsch is
sufficient to guarantee a rooted $H$-minor.
If so, we say that $H$ is {\em \rcn}.
(We say that $H$ is {\em weakly \rcn} if an \Shsch guarantees only an $H$-minor.)

In connection with Hadwiger's conjecture Duchet and Meyniel~\cite{D,DM} studied a related notion of scheme, which we will call a Duchet-scheme to avoid confusion:
An {\em \hDschn} in $G$ is a set of
connected subgraphs $\{C_e: e\in E(H)\}$ of $G$ such that $C_e\cap C_f\neq\emptyset$
if and only if $e$ and $f$ have a common endpoint (in $H$).  They conjectured
that every $p$-chromatic graph contains a $K_p$-\Dschn.  They observed that $G$
contains an \hDsch if and only if $G$ contains a {\em simple \hDschn}, which is
defined like an \Shsch except that Definition~\ref{d:scheme-def}.2 is
required only for $k=2$.  

Observe that $K_{1,3}$ contains a simple $K_3$-\Dsch but no $K_3$-minor.
In general a \Dsch does not guarantee a minor for any graph with a $K_3$-subgraph,
leading us to the strengthening in Definition~\ref{d:scheme-def}.2.
\medskip

In the next section, we collect basic definitions from graph theory and introduce terminology we use in the paper. We also establish some observations about \Shschs that are used to introduce the more restrictive ``\Chschn'' in Section~\ref{forest-sec} that is equivalent for our purposes.
 We prove that many graphs that are \rcn, including
 forests (Section~\ref{forest-sec}),
 cycles (Section~\ref{cycle-sec}),
 subgraphs of $K_4$ and $K_{1,1,3}$ (Sections~\ref{complete-sec} and~\ref{K23-sec}),
 and certain graphs built from such blocks;
 Corollary~\ref{summary} summarizes all these positive results.
 On the other hand, we develop a general technique
 in Section~\ref{short-sec} that allows us to show that many graphs are not \rcn,
 including several 7-vertex graphs (discussion following Theorem~\ref{theta}),
 certain graphs with few cycles (Corollary~\ref{2odds} and Theorem~\ref{theta}),
 and graphs that are not 6-colorable (Theorem~\ref{chromatic bound}).
 The results for $K_7$ and cycles were also obtained independently by Chan and Seymour~\cite{Se}, and the content of Section~6 by Seymour and Thomas~\cite{Se}, but these efforts were never formally written down.  We also show that $K_5$ and $K_{3,3}$ are weakly \rc (Theorem~\ref{k5-k33}) and we construct graphs that are not weakly \rc (following Theorem~\ref{theta}).
 In Section~\ref{question-sec}, we ask some open questions.

 We start by clarifying some basic notions.

\section{Definitions and observations}\label{def-sec}

For basic definitions from graph theory, we refer the reader to the introductory graph theory book by Diestel~\cite{Di}.

All graphs in this paper are simple, that is we do not allow loops or multiple edges.
A {\em $u,v$-path} is a connected graph in which $u,v$ have degree 1 and all other vertices (the {\em internal} vertices) have degree 2. Equivalently we may refer to a $u,v$-path as a sequence of distinct vertices starting at $u$ and ending at $v$ in which consecutive vertices are adjacent.

We say that $G$ contains an {\em $H$-immersion} if $H$ can be obtained from $G$ by a sequence of vertex-deletions, edge-deletions, and edge-liftings (that is, replacing the edges in a path $uvw$ by an edge $uw$). It is not difficult to see that  if $G$ contains an \Shsch in which all paths are edge-disjoint, then $G$ contains an $H$-immersion.
(The converse is false, since $K_{1,4}$ contains a $2K_2$-immersion, but no $2K_2$-scheme.)

When we {\em contract} an edge $uv$ in a graph $G$
we obtain a new graph with vertex set $V(G)\cup\{w\}\setminus\{u,v\}$ (where $w$ is a new
vertex) and the same edges as $G$, except that $wx$ is an edge if at least one of $ux$ or $vx$ is an edge in $G$.

We say that a graph $H$ is a {\em minor} of another graph $G$, or that $G$ {\em contains an $H$-minor}, if a graph isomorphic to $H$ can be obtained from $G$ by a sequence of vertex-deletions, edge-deletions and edge-contractions. Equivalently, there is a collection of vertex-disjoint connected subgraphs $C_v$ in $G$ (one for each $v\in V(H)$), such that for each edge $uv$ in $H$ there is
an edge connecting $C_u$ and $C_v$. If $V(H)\subseteq V(G)$
and for all $v\in V(H)$ we have $v\in V(C_v)$, then we say that $G$ has a {\em rooted} $H$-minor.

It is obvious that if $H$ is a rooted minor of $G$, then we can
find an \Shsch in $G$. We are concerned with cases in which
the converse holds: when can the paths in the \Shsch be
``rerouted'' or ``untangled'' to get a (rooted) $H$-minor?

Since we are interested in whether the presence of an \Shsch is
sufficient to guarantee a (rooted) $H$-minor,
we henceforth assume that $G$ is the underlying graph of an \Shschn.

\begin{definition}
We call a graph $H$ {\em \rcn} if every graph
$G$ that contains an \Shsch must have a rooted $H$-minor.
We call a graph $H$ {\em weakly \rcn} if every graph
$G$ that contains an \Shsch must have an $H$-minor.
\end{definition}

Every \rc graph is weakly \rcn.
$K_1$ and $K_2$ are trivially \rcn, and we will see that many other
graphs are \rcn.

The following
observation implies that there is a forbidden subgraph characterization for
\rc graphs.

\begin{lemma}\label{subgraph}
Every subgraph of a \rc graph is \rcn.
\end{lemma}

\begin{proof}
Let $H$ be a \rc graph and let $H'$ be a subgraph of $H$.
Suppose that there is an \Shpsch in a graph $G'$.  We must show
that $H'$ is a rooted minor of $G'$.
Define $G$ to be the graph $G$ with
vertex-set $V(G)=V(G')\cup (V(H)\setminus V(H'))$ 
and edge-set $E(G)=E(G')\cup (E(H)\setminus E(H'))$.
We can extend the \Shpsch in $G'$ to an \Shsch in $G$ by letting
$P_{uv}=uv$ for every edge $uv\in E(H)\setminus E(H')$. Now, $G$ has a rooted
$H$-minor, and each edge $uv\in E(G)\setminus E(G')$ cannot be involved in it
except as the edge between $C_u$ and $C_v$.  Also, for each vertex $v\in
V(G)\setminus V(G')$, clearly $C_v=\{v\}$.  Therefore, restricting the rooted
$H$-minor in $G$ to the subgraph $G'$ yields a rooted $H'$-minor in $G'$.
\end{proof}

Unfortunately the proof of this basic result does not extend to the weakly
\rc case, and so we primarily focus our attention on the \rc case.

The next simple observation shows that minimal non-\rc graphs are connected,
which allows us to restrict our attention to the case when $H$ and $G$ are connected.
Note that Lemma~\ref{components} cannot be strengthened to 2-connected components,
since the 9-vertex graph consisting of two 5-cycles with one shared  vertex
is an example of a graph that is not contractible (by Corollary~\ref{2odds}) even though
its blocks are contractible (by Theorem~\ref{cycle}).

\begin{lemma}\label{components}
A graph is (weakly) \rc if and only if every one of its components
is (weakly) \rcn.
\end{lemma}

\begin{proof}
The forward implication is obvious in the \rc case, since every
component is a subgraph. For the forward implication of the weakly
\rc case, suppose that $H$ is weakly \rcn, $H'$ is a
component of $H$, and $G'$ is a graph that contains an \Shpschn.  Let $G$
be obtained from $H$ by replacing the component $H'$ by $G'$ (or equivalently,
construct $G$ as in the proof of Lemma~\ref{subgraph}).  Then $G$ contains
an \Shschn, where edges of $H'$ are represented by paths in $G$ from the
\Shpsch in $G'$, and edges in $E(H)\setminus E(H')$ are represented by
themselves. 
Since $H$ is weakly \rcn, for each
component $C$ of $H$ there is a $C$-minor $M_C$ in $G$, and these minors
in $G$ are pairwise disjoint.  Define an auxiliary digraph $D$
with vertex set $\{$components of $H\}\cup\{G'\}$
($=\{$components of $G\}\cup\{H'\}$), such that for each
component $C$ of $H$ we have an edge from $C$ to the component $C'$ of $G$
that contains $M_C$.  In $D$, every component of $H$ has out-degree 1, $G'$ has out-degree 0 and $H'$ has in-degree 0.
Observe that the vertices on a cycle in $D$ must be isomorphic components of $H$ and have in-degree 1 in $D$.
Thus if we start a walk at $H'$ in $D$, then it cannot contain a cycle in $D$ (the first vertex on the walk that is in this
cycle would not have in-degree 1 in $D$), and must therefore end in $G'$. Since the minor-relation is transitive,
this path from $H'$ to $G'$ in $D$ implies that $G'$ contains a minor of $H'$, as desired.

For the reverse implications, let $H$ be a graph, all of whose components are
(weakly) \rcn. Consider an \Shsch in $G$.  For each
component $C$ of $H$, let $G_C$ be its underlying graph, i.e. the union of the scheme-paths $P_e$
with $e\in E(C)$.  Clearly $G_C$ contains a $C$-scheme, and since $C$ is
(weakly) \rcn, $G_C$ contains a ($C$-minor or) rooted $C$-minor.
By the definition of an \Shschn, for distinct components $C,C'$ of $H$, $G_C$
and $G_{C'}$ are vertex-disjoint.  Thus the subgraphs $G_C$ for all components $C$
of $H$ can be combined to get (an $H$-minor or) a rooted $H$-minor in $G$.
\end{proof}

Note that by Definition~\ref{d:scheme-def}.2, for each
$w\in V(G)\setminus V(H)$, there is a vertex $f(w)=u\in V(H)$ such that $w$
only intersects the \Shsch in paths starting at $u$.  Moreover,
Definition~\ref{d:scheme-def} could be restated in an equivalent form if
we replaced Condition~\ref{d:scheme-def}.2 by by the statement that there exists a function
$f: V(G)\to V(H)$ such that for each $v\in V(H)$ we have $f(v)=v$ and
$f(V(P_{uv}))=\{u,v\}$.  
We continue our investigation in this direction in the next section.


\section{Colored schemes and forests}\label{forest-sec}

We will see that
whenever $G$ has an \Shschn, there is a rooted minor of $G$ that contains a
``\Chschn'', which is an \Shsch with additional properties.  It will follow that, for the
purpose of investigating whether a graph $H$ is (weakly) \rcn, we may
restrict our attention to \Chschn s.

\begin{definition}\label{d:colored-def}
A  {\em \Chschn}
is an \Shsch with underlying graph $G$ and a map $f:V(G)\to V(H)$ such that
\begin{enumerate}
\item $f(v)=v$ for each $v\in V(H)$,\label{cd1}
\item $f(V(P_{uv}))=\{u,v\}$ for each $uv\in E(H)$,\label{cd2}
\item $f$ is a proper coloring of $G$, and\label{cd3}
\item each vertex in $V(G)\setminus V(H)$ has degree at least 4.\label{4deg-def}
\end{enumerate}
\end{definition}

\begin{remark}\label{r:colored-equiv}{\rm
We collect some obvious properties of \Cschn s:

\begin{enumerate}
\item The colors on every path $P_{uv}$ alternate between $u$ and $v$.
\item Every edge in $G$ is in exactly one $P_{uv}$.
\item Every $P_{uv}$ is an induced subgraph of $G$.  (This follows from the
previous two.)
\item $f$ is a homomorphism from $G$ to $H$ that fixes $V(H)$.
\item The scheme-paths form an $H$-immersion in $G$.
\item Every vertex in $V(G)\setminus V(H)$ is in at least 2 scheme-paths.\label{4deg-rem}
\item If $u$ has degree 2 in $H$, then $P_{uv}$ contains every vertex of color $u$.\label{2deg-rem}
\end{enumerate}
}\end{remark}

Thus, a \Chsch is
a vertex-coloring of $G$ such that the
vertices of $H$ receive different colors, and each path $P_{uv}$ is a
Kempe chain connecting $u$ and $v$.\footnote{Motivated directly by Kempe
chains, Duchet and Meyniel~\cite{D,DM}
 define a ``chromatic $p$-scheme'' to be a simple
$K_p$-\Dsch that satisfies Conditions~\ref{d:colored-def}.\ref{cd1},
\ref{d:colored-def}.\ref{cd2}, and \ref{d:colored-def}.\ref{cd3}.
They showed that a $p$-chromatic graph does not necessarily have
a chromatic $p$-scheme, making this definition more restrictive than their
``$K_p$-\Dschn'' and ``simple $K_p$-\Dschn'' for their interest in
Hadwiger-like conjectures.}
Remark~\ref{r:colored-equiv}.\ref{4deg-rem} ensures that these Kempe chains contain no topologically irrelevant vertices.

The following lemma shows that for the purpose of this paper, \Sschs and \Cschs are equivalent notions.

\begin{lemma}\label{colored-lem}
Let $G,H$ be graphs with $V(H)\subseteq V(G)$.
$G$ has an \Shsch if and only if some rooted minor of $G$
is the underlying graph of a \Chschn.
\end{lemma}

\begin{proof}
For the backward implication observe that if $H$ has a \Csch in some
rooted minor $M$ of $G$, then we can find an \Shsch in $G$ by reversing all contractions
performed to obtain $M$ from $G$, potentially lengthening the paths of the \Shsch
in $M$. For the forward direction, let $R=\{P_e:e\in E(H)\}$ be an \Shsch in $G$.
We convert $R$ into a \Csch of $H$ in some rooted minor of $G$, by repeatedly
performing the following operations, with priority given to operations listed earlier.

\begin{enumerate}
\item If $G$ contains vertices or edges that are not used in $R$,
then we delete them to ensure that $G$ is the underlying graph of $R$.
\item If $G$ contains a vertex that is not in $V(H)$ and
that is on only one path in $R$, then we can contract one of
its incident edges. This will shorten only this one path, but leave
all other paths unaffected.
\item If any $P_{uv}\in R$ is not an induced path subgraph of $G$, then we can
replace $P_{uv}$ in $R$ by a proper subsequence of its vertices that yields an induced $uv$-path. (This may
enable us to delete or contract additional edges or vertices.)
\item By Operations 1 and 2 and Definition~\ref{d:scheme-def}.2, there is a unique function
$f:V(G)\to V(H)$ such that $f(v)=v$ for $v\in V(H)$, and
for each $v\notin V(H)$, $f(v)$ is the
common endpoint of two paths in $R$ that contain $v$. If there is an edge $uv$
with $f(u)=f(v)$, then we contract $uv$. This will shorten some paths and
may create chords on others, in which case we repeat the previous steps.
\end{enumerate}

This process terminates with a minor $M$ of $G$, since every time we execute
Operation~1, 2, or 4, we shrink $G$, and Operation~3 decreases the size of a
path in $R$.
It is straightforward to check that when we execute Operation~4 for the last time,
the $f$ we obtain satisfies all conditions in Definition~\ref{d:colored-def} for $M$.
\end{proof}

In order to prove that a graph $H$ is \rc (weakly \rcn),
by Lemma~\ref{colored-lem} it suffices
to show that given any \Chschn, its underlying graph $G$ contains
a rooted $H$-minor (an $H$-minor).

\begin{theorem}\label{leaf}
If $H$ has a vertex $u$ of degree 1, then $H$ is \rc if and
only if $H\graphminus u$ is \rcn.
\end{theorem}

\begin{proof} The forward implication follows from Lemma~\ref{subgraph}.
For the reverse implication, suppose that $H\graphminus u$ is \rc and
consider a \Chsch in $G$; by Lemma~\ref{colored-lem} it suffices
to prove that $G$ contains a rooted $H$-minor.
Let $v$ be the unique neighbor of $u$ in $H$.
By Remark~\ref{r:colored-equiv}.\ref{4deg-rem}, the only vertex of color $u$ in $G$
is $u$ itself; therefore, $P_{uv}$ has only one edge. Thus,
$G\graphminus u$ has a \Csch of $H\graphminus u$, so it contains
$H\graphminus u$ as a rooted minor. Adding the edge $uv$ in $G$ to it,
we obtain a rooted $H$-minor in $G$.
\end{proof}

It follows by induction that all trees are \rcn.
Thus, by Lemma~\ref{components},
we obtain:

\begin{cor}
Forests are \rcn. \qed
\end{cor}


\section{Contracting cycles}\label{cycle-sec}

The next result is similar in flavor to Theorem~\ref{leaf}.

\begin{lemma}\label{triangle}
 If $H$ contains a triangle $uvw$ and
$u,v$ each have degree 2 in $H$, then $H$ is \rc if and only if
$H\graphminus\{u,v\}$ is \rcn.
\end{lemma}

\begin{proof}
Again one implication is obvious, so we will now assume
that $H'=H\graphminus\{u,v\}$ is \rcn. Consider a \Chsch
in $G$. Obtain $G'$ by removing all vertices of
colors $u$ and $v$ from $G$. Then $G'$ contains a \Csch of
$H'$, so $G'$ contains a rooted $H'$-minor. To extend this $H'$-minor to
a rooted $H$-minor of $G$, observe that $w$ is adjacent to a vertex $u'$ of color $u$,
and a vertex $v'$ of color $v$. By Remark~\ref{r:colored-equiv}.\ref{2deg-rem}, $u'$ and $v'$  must
be on $P_{uv}$ and
we can contract $P_{uv}$ to the edge $u'v'$ to obtain
a rooted $H$-minor in $G$.
\end{proof}

Lemma~\ref{triangle} implies that triangles are \rcn. With a bit more effort we obtain
the next result which has also been obtained independently by Melody Chan and Paul
Seymour~\cite{Se}. Our proof uses the same idea that was used in~\cite{FKPR} to show
that 4-cycles are \rcn.

\begin{theorem} \label{cycle}
Cycles are \rcn.
\end{theorem}

\begin{proof}
Let $n\geq 3$ and let $C_n$ be a cycle with vertices $u_1, u_2, \dots,
u_n$ (in cyclic order). We will show by induction on $|V(G)|$ that if a
graph $G$ contains a \Csch of $C_n$, then $G$ contains $C_n$ as a rooted minor.
If $G$ has $n$ vertices, then $G$ contains $C_n$ as a subgraph, which suffices.

Now suppose that $G$ has more than $n$ vertices.
By Lemma~\ref{colored-lem}, we can assume that $G$ is the underlying graph of a colored $C_n$-\Sschn.
By Remark~\ref{r:colored-equiv}.\ref{2deg-rem}, each path
$P_{u_iu_{i+1}}$ contains all vertices of colors $u_i$ and $u_{i+1}$.
Thus there are equally many vertices of color $i$ and $i+1$, and
it follows that all colors appear equally often. Then, because
$G$ has more than $n$ vertices, it follows that every $P_{u_iu_{i+1}}$ is
non-trivial. Thus for each $i$, the neighbor of $u_i$ on $P_{u_iu_{i+1}}$ (call it
$v_i$) is different from $u_{i+1}$.
Contract all edges of the form $u_iv_i$ in $G$ to obtain a new graph $G'$.
Since $P_{u_iu_{i+1}}$ contains all vertices of color $u_i$, $P_{u_iu_{i+1}}$
contains $v_{i-1}$. Hence, after contraction, $P_{u_iu_{i+1}}$ contains a path
$P_{u_{i-1}u_{i}}'$ from $u_{i-1}$ to $u_i$.  If we recolor each vertex of color
$u_i$ in $V(G')\setminus V(C)$ by color $u_{i-1}$, for all $i$, then each internal
vertex of $P_{u_{i-1}u_{i}}'$ has color $u_{i-1}$ or $u_i$; hence, these paths
form a $C_n$-scheme in $G'$.  By induction,
$G'$ contains a rooted $C_n$-minor.  Since $G'$ was obtained from $G$ by
contractions, $G$ contains a rooted $C_n$-minor as well.
(Uncontract to obtain it.)
\end{proof}

A {\em block} in a graph is a maximal subgraph that has no cut-vertex.
A {\em cactus} is a connected graph in which every block is an edge or a cycle.
Combining Theorem~\ref{cycle}, Lemma~\ref{triangle}, and Theorem~\ref{leaf}, we obtain
a small extension of the fact that trees are \rcn.

\begin{cor}\label{cactus}
A cactus in which at most one cycle has more than 3 vertices is \rcn. \qed
\end{cor}

\section{Small cliques and bicliques}\label{complete-sec}

We have already seen that $K_n$ is \rc for $n\le 3$.
We will extend this to $n=4$, using
the characterization of graphs without rooted $K_4$-minors given by Robertson, Seymour and Thomas~\cite{RST}.
The results in this section were also obtained independently by Robin Thomas and Paul
Seymour~\cite{Se}.

\begin{theorem}\label{K4}
$K_4$ is \rcn.
\end{theorem}

\begin{proof}
If $K_4$ is not \rcn, consider a minimal colored $K_4$-scheme which has no rooted $K_4$-minor.  Let $Z=\{z_1,z_2,z_3,z_4\}$ be the set of roots and let $G$ be the underlying graph. According to (2.6), p288 in~\cite{RST}, either (i) $G$ has a rooted $K_4$-minor, or
(ii) there are sets $A_1,A_2,B$ such that $G[A_1]\cup G[A_2]\cup G[B]=G$, $A_1\cap A_2=A_1\cap B=A_2\cap B$
with $|A_1\cap A_2|=2$, and, without loss of generality, $Z\cap (A_i\setminus B)=\{z_i\}$ for $i=1,2$, or
(iii) there are sets $A,B$ such that $G[A]\cup G[B]=G$, $|A\cap B|\le 3$,
$Z\subseteq A$, $|B\setminus A|\ge 2$, and $|Z\cap B|\le 2$, or
(iv) $G$ can be drawn in a plane so that every vertex of $Z$ is incident with the infinite region.

Case (i) is an immediate contradiction.  In case (iv), if $z_1,z_2,z_3,z_4$ occur in this order
on the boundary of the infinite region, then $P_{z_1z_3}$ and $P_{z_2z_4}$ must cross, another contradiction.

In case (ii), $A_1\cap A_2$ must contain vertices $z_1',z_2'$ of colors~$z_1$ and $z_2$, respectively, for use by the paths in the colored $K_4$-scheme.
Let $Z'=\{z_1',z_2',z_3,z_4\}$.
Restricting the colored $K_4$-scheme to $G[B]$ gives a new colored $K_4$-scheme with root set $Z'$,
and by minimality this has a rooted $K_4$-minor $H$.
Restricting $P_{z_1z_3}$ and $P_{z_2z_4}$ from the original colored $K_4$-scheme
to $A_1$ and $A_2$, and adding this to $H$ we obtain a rooted minor of $K_4$ in $G$ with
roots in $Z$.

Thus we may assume that we have case (iii).
Let $v$ be any vertex in $B\setminus A$. Each scheme-path through $v$ intersects
$A\cap B$ in two vertices, so $|A\cap B|\ge 2$.  Moreover, $v$ is in at least
two such paths and $|A\cap B|\le 3$, so two such paths must have a common
vertex $x\in A\cap B$ of the same color as $v$.
If $A\cap B=\{x,y\}$, then $G[A]\cup\{xy\}$ contains a $K_4$-scheme
with the same roots, and thus by minimality $G[A]\cup\{xy\}$ contains a rooted $K_4$-minor,
which we can extend to
the desired $K_4$-minor in $G$ by replacing $xy$ with the $x,y$-portion of a path in the
\Csch that goes through $v$. So we may assume that $A\cap B=\{x,y,z\}$. Since $v$
is not adjacent to $x$ and $d_G(v)\ge 4$, it follows that $v$ has at least two
neighbors in $G\graphminus A$. This argument is independent of the choice of $v\in B\setminus A$, so $G\graphminus A$
contains a cycle $C$.
There must be 3 disjoint paths from $\{x,y,z\}$ to $C$,
since otherwise there would be sets $A',B'$ such that $A\subseteq A'\setminus B'$, $V(C)\subseteq B'\setminus A'$,
$G[A']\cup G[B']=G$, and $|A'\cap B'|\le 2$, and we just showed how to handle this case.
Now let $G'=G[A]\cup\{xy,yz,xz\}$ and observe that again $G'$ has a $K_4$-scheme with
roots $Z$, so by minimality $G'$ contains a rooted $K_4$-minor $H$.  We can
replace $xy,yz,xz$ in $H$ by $C$ and the 3 disjoint paths from $\{x,y,z\}$ to
$C$ as necessary, to obtain a rooted $K_4$-minor in $G$.
\end{proof}

The next lemma is useful for studying the case $n=5$.  A {\em pasting (along $K_m$)} of two graphs $G_1,G_2$ is the graph
$G_1\cup G_2$ if $V(G_1)\cap V(G_2)$ induces a complete graph $K_m$ in both $G_1$ and $G_2$.

\begin{lemma}\label{pasting}
Let $H$ be a graph with no cut-set with at most $m$ vertices.

\begin{enumerate}
\item If the pasting of two graphs $G_1,G_2$ along $K_m$ contains a \Chschn,
then $G_1$ or $G_2$ contains a \Ssch of a graph isomorphic to $H$.
\item If every edge-maximal graph with no $H$-minor is obtained by repeated pastings
along $K_k$ with $k\le m$ of graphs
that do not contain \Cschn s of graphs isomorphic to $H$, then $H$ is weakly \rcn.
\end{enumerate}
\end{lemma}

\begin{proof} For the first statement,
let $S$ be the set of vertices in $H$ for which the color $v$ appears in
$V(G_1)\cap V(G_2)$.  Then $|S|\le |V(G_1)\cap V(G_2)|=m$, so $S$ is not a
cut-set of $H$.  For any $u,v\in V(H)\setminus S$, $H\graphminus S$ contains a
$u,v$-path $Q$.  The union of scheme-paths $\bigcup\{P_e: e\in E(Q)\}$ is a
connected subgraph of $G$ that contains $u$ and $v$ but no vertices with
colors in $S$, so it is all in one component of $(G_1\cup G_2) \graphminus
(V(G_1)\cap V(G_2))$.  Thus, without loss of generality we may assume that
$V(G_1)\setminus V(G_2)$ contains 
every vertex of the \Chsch whose color is not in $S$.

Replace each vertex $u\in V(H)\setminus V(G_1)$ with any vertex $u'$ of the same
color in $V(G_1)\cap V(G_2)$.  For each scheme-path $P$ ending at $u$,
replace $P\graphminus V(G_1)$ by a single edge to $u'$, which can be done since
$V(G_1)\cap V(G_2)\cong K_m$.  Similarly, we can replace
every other intersection of a scheme-path with $V(G_2)\setminus V(G_1)$ by a single edge.
The result contains an \Shpschn, where $H'\cong H$ and the vertices $u$ in $V(H)\setminus V(G_1)$
have been replaced with vertices $u'$ in $V(G_1)\cap V(G_2)$.

To prove the second statement, suppose that $H$ is not weakly \rcn. So there exists a graph
that contains an \Shsch but no $H$-minor; let $G$
be such a graph with the minimum number of vertices. By Lemma~\ref{colored-lem} we
can assume that $G$ is the underlying graph of a \Chschn.  Extend $G$ to an edge-maximal
$H$-minor free graph $G'$. By Lemma~\ref{pasting}.1, the vertex-minimality of $G$
implies that $G'$ is not a pasting of two smaller graphs. But this contradicts our hypothesis
and the fact that $G$ has a \Chschn.
\end{proof}

The Hanani-Tutte theorem~\cite{Hanani,Tutte} (see~\cite{PSS} for a history and short proof) says that if a graph
is drawn in the plane such that every pair of non-adjacent edges crosses an even number of times,
then the graph is in fact planar.  If the underlying graph
of a \Chsch is planar, then this corresponds to a drawing of $H$ in the plane in which pairs
of non-adjacent edges in $H$ do not intersect, so by the Hanani-Tutte theorem, $H$ must be planar.
Thus there are no planar \Cschs for $K_5$ or $K_{3,3}$, and it follows that a graph is planar
if and only if it does not contain a \Csch of $K_5$ or $K_{3,3}$.  (Duchet and Meyniel observed this
for (simple) \Dschs by the same arguments~\cite{D,DM}.)
It follows that a \Csch of $K_5$ or $K_{3,3}$ must contain a minor of either $K_5$ or $K_{3,3}$, which indicates that $K_5$ and $K_{3,3}$ are likely to be \rcn. In fact, we easily obtain

\begin{theorem}\label{k5-k33}
$K_5$ and $K_{3,3}$ are weakly \rcn.
\end{theorem}

\begin{proof}
By a theorem of Wagner~\cite{W} (see Theorem 8.3.4 in Diestel~\cite{Di}), a maximal graph without a $K_5$-minor
can be obtained by pasting planar triangulations and
copies of the 3-regular M\"obius ladder $M_8$ (obtained by pairwise making opposite vertices of an
 8-cycle adjacent) along triangles and edges. Since a colored $K_5$-scheme has minimum degree 4,
 $M_8$ cannot contain this, and the result follows by Lemma~\ref{pasting}.2.

 The proof for $K_{3,3}$ follows from Wagner's result that every edge-maximal graph without a $K_{3,3}$-minor can be constructed from pasting copies of $K_5$ and planar triangulations along edges.
\end{proof}

Unfortunately, there is no corresponding characterization of edge-maximal $K_6$-minor free
graphs. The topological concept that corresponds to the non-planarity of $K_5$ is the
fact that $K_6$ has no linkless embedding in 3-space~\cite{CG,S}. Using a similar approach as above
it can be shown there is no colored $K_6$-scheme that can be linklessly embedded in 3-space~\cite{Se}, a strong
indication that $K_6$ is weakly \rcn, and possibly even \rcn. In the next section
we will show that $K_7$ is not \rcn, but we have no indication that $K_n$ is not
weakly \rc even for much larger values of $n$. It remains open whether $K_{2,4}$ or $K_{3,3}$
are \rcn, and we are unaware of any bipartite graphs that are not \rcn.

\section{A contraction lemma}\label{K23-sec}

We now extend the idea from Theorem~\ref{leaf} to obtain a useful technical lemma.
\smallskip

\begin{lemma}\label{removable}
Let $H$ be a graph with a stable set $S$ and a forest subgraph $F$ with $V(F)=S\cup N(S)$,
such that isolated vertices of $F$ are in $S$, and every other component of $F$ is a star
with a root in $N(S)$ and non-root leaves in $S$.
Suppose $G$ contains a \Chschn. If
\begin{enumerate}
\item $H\graphminus S$ is \rcn, and
\item for each $u\in S$ and $w\in N(u)$ with $uw\not\in E(F)$,
the second vertex on $P_{uw}$ is on
$P_{vw}$ for some $vw\in E(F)$,
\end{enumerate}
then $H$ is a rooted minor of $G$.
\end{lemma}

\begin{figure}[ht]
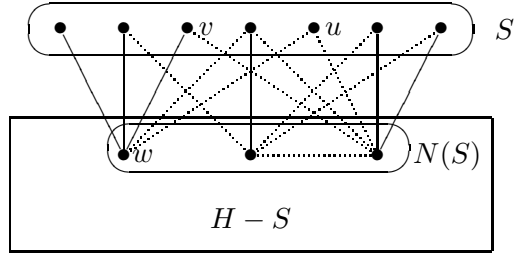

\centerline{
\hbox{\xy /r2pc/:,
(0,6)*{\bullet}="S1", (1,6)*{\bullet}="S2", (2,6)*{\bullet}="S3", (3,6)*{\bullet}="S4",
(4,6)*{\bullet}="S5", (5,6)*{\bullet}="S6", (6,6)*{\bullet}="S7", (7,6)*{S},
(1,4)*{\bullet}="T1", (3,4)*{\bullet}="T2", (5,4)*{\bullet}="T3", (6.1,4)*{N(S)},
(2.3,6)*{v},(4.3,6)*{u},(1.3,4)*{w},
"S1";"T1"**@{-},"S2";"T1"**@{-},"S3";"T1"**@{-},
"S4";"T2"**@{-},"S6";"T3"**@{-},"S7";"T3"**@{-},
"S4";"T1"**@{.},"S2";"T2"**@{.},"S3";"T3"**@{.},"S5";"T1"**@{.},
"S4";"T3"**@{.},"S5";"T3"**@{.},"S7";"T2"**@{.},"S6";"T2"**@{.},
"T3";"T2"**@{.},
"S1"+<-1pc,-0.8pc>;"S7"+<1pc,0.8pc> **\frm<44pt>{-},
"T1"+<-0.5pc,-0.5pc>;"T3"+<1pc,1pc> **\frm<44pt>{-},
(-0.8,4.6); (6.8,2.5) **\frm{-}, (3,3) *{H\graphminus S}
\endxy}}
\caption{\label{removable pic} The stable sets $S$ in  Lemma~\ref{removable}, solid edges are in $F$. }
\end{figure}

\begin{proof}
For each edge $vw$ in $F$ with $v\in S$, contract $P_{vw}\graphminus v$ to $w$.
Then remove all remaining vertices in $V(G)\setminus S$ that have colors in $S$.
Let $G'$ be the resulting graph.

For each $w\in N(S)$, let $V_w=\bigcup_{vw\in E(F)}(V(P_{vw})\setminus v)$, the set
of vertices of $G$ that are contracted to $w$.  The colors in $V_w$ lie in $\{w\}\cup N_F(w)$,
so for distinct $w_1,w_2$ in $N(S)$, the sets $V_{w_1},V_{w_2}$ do not intersect.
Therefore, the contraction step of the construction is well-defined.  Also,
$V_w\cap V(H)=\{w\}$, so vertices of $V(H)$ remain distinct after contraction.

After the contractions, each scheme-path $P_{xy}$ will still be connected, so
it will contain an $x,y$-path $P_{xy}'$ with colors in $\{x,y\}$.
For each $xy\in E(H\setminus S)$, $G'\graphminus S$ contains $P_{xy}'$, so
$G'\graphminus S$ contains a colored $(H\graphminus S)$-scheme.  Since $H\graphminus S$ is \rcn, $G'\graphminus S$
contains a rooted $H\graphminus S$-minor.

For $u\in S$ and $uw\in E(H)$, the second vertex of $P_{uw}$ is contracted to $w$,
by construction if $uw\in E(F)$ and by Condition~\ref{removable}.2
if $uw\not\in E(F)$.  Therefore $P_{uw}'$ is a single-edge path, which is contained in $G'$.
Hence, for $u\in S$ and $uw\in E(H)$, we have $uw\in E(G')$, and
$u$ and $uw$ are not in $G'\graphminus S$.
Therefore we can add $S$ and incident edges in $H$ to the rooted $H\graphminus S$-minor in $G'\graphminus S$,
and obtain a rooted $H$-minor in $G'$.  Since $G'$ was obtained from $G$ by contractions
and deletions, $G$ must contain a rooted $H$-minor as well.
\end{proof}

This lemma easily implies Theorem~\ref{leaf} by letting $S=\{u\}$ and $F=\{uv\}$ (where $u$ is a leaf with neighbor $v$) and observing that
Condition~\ref{removable}.2 is vacuously satisfied. Similarly it follows that even cycles are \rcn, by letting
$S$ be one of the partite sets and $F$ be a perfect matching  and observing that Condition~\ref{removable}.2 
follows from Remark~\ref{r:colored-equiv}.\ref{2deg-rem}. With more care, we obtain:

\begin{theorem}
$K_{2,3}$ and $K_{1,1,3}$ are \rcn.
\end{theorem}

\begin{proof}
It suffices to consider $H=K_{1,1,3}$, since it contains $K_{2,3}$ as a subgraph.
Let $V(H)=\{a,b,c,x,y\}$, where $T=\{a,b,c\}$ is the stable set of size 3. Consider a \Chsch
with underlying graph $G$.

If there is a scheme-path $P_{uv}$ with $u\in T$ that has only one edge, then by
Remark~\ref{r:colored-equiv}.\ref{2deg-rem} there is only one vertex of color $u$ and hence $ux,uy\in E(G)$.
Thus $G\graphminus u$ contains a \Csch of $H\graphminus u$. Since $H\graphminus u$ is a subgraph of $K_4$ and thus \rcn,
$G\graphminus u$ contains $H\graphminus u$ as a rooted minor, so that $G$ has the desired $H$-minor.

Thus for each scheme-path $P_{uv}$ with $uv\neq xy$, we can
let $v_u$ denote the second vertex on $P_{uv}$; note that the color of $v_u$ is $v$, but $v_u\neq v$.
By Remark~\ref{r:colored-equiv}.\ref{4deg-rem}, $v_u$ is on at least one more scheme-path,
which has the form $P_{*v}$.

Suppose without loss of generality that $x_a$ is also on $P_{bx}$.
If $y_b\in P_{ay}$, then we can apply Lemma~\ref{removable}
with $S=\{a,b\}$ and $F=\{ay,bx\}$ since $H\graphminus S$ is a triangle (and thus \rcn), and the result follows. So we may assume that $y_b\in P_{cy}$. Next, we may assume that $x_c\in P_{ax}$, since the possibility $x_c\in P_{bx}$ reduces similarly with $S=\{b,c\}, F=\{bx,cy\}$.
Continuing like this, it follows that $y_a\in P_{by}$, $x_b\in P_{cx}$, and $y_c\in P_{ay}$.

Recolor the vertices of $V(G)\setminus V(H)$ by changing color $b$ to $a$, color $c$ to $b$,
and color $a$ to $c$. Since $x_a\in V(P_{bx})$ we can let $P_{ax}'$ be the path formed by $a$, $a x_a$, and the
$x_a,x$-subpath of $P_{bx}$; this and similarly-defined paths form a \Chsch
for the new coloring of $G$.  Since each edge in the $x_a,x_b$-path of $P_{bx}$ has
endpoints of colors $a$ and $x$ in the recoloring, these edges are not in any path
of the new \Chschn, so we can remove the edges and obtain the result by induction.
There are no such edges if $x_a=x_b$, but then $x_b$ is in $P_{ax}$, a case that we
ruled out earlier.
\end{proof}

\section{Colored schemes with short paths}\label{short-sec}

We can summarize our knowledge about \rc graphs as follows.

\begin{cor}\label{summary}
Suppose that $H$ is a graph where every block is a $K_3$, $K_2$, or $K_1$, except
for up to one block per component, which can be a $K_4$, $K_{1,1,2}$, $K_{1,1,3}$, $K_{2,3}$, or a cycle (of any length).
Then $H$ is \rcn.
\end{cor}

To find graphs that are not \rcn, it seems unclear if one should
consider \Sschs in which the paths are short and offer little
choice for contraction, or long and entangled in
complicated ways. In specific examples we found that as
soon as the paths get long one can typically find the desired minor,
even though it may be hard to prove this in general. Thus we focus on the case when
all paths in the \Ssch are length 3 (short, but non-trivial.) Specifically,
we consider the case when there is only one additional vertex in each
color, and $M'(H)$ will denote the underlying graph.

\begin{definition} Let $H$ be a graph. Let $M'(H)$ be the
rooted graph obtained as follows: for every $v\in V(H)$, $M'(H)$
contains two vertices $v_1,v_2$ and $u_i$ is adjacent to
$v_j$ whenever $u$ is adjacent to $v$ in $H$ and $i=2$ or
$j=2$ (or both).

We say that $H$ is {\em \mpc} if $M'(H)$ contains a rooted
$H$-minor, where $v_1$ is the root in $M'(H)$ corresponding to $v\in V(H)$.
\end{definition}

The notation $M'(H)$ is derived from Mycielski's construction;
except for its ``center vertex'', this is the same construction.
Note that $H$ is trivially ``weakly \mpc'',
since $M'(H)\graphminus V(H)$ is isomorphic to $H$.

Observe that $M'(H)$ contains a \Csch of $H$
where $v=v_1$ and the color of
$v_2$ is $v$: for every edge $uv$ in $H$, there is the
corresponding path $u_1v_2u_2v_1$ in $M'(H)$.
Thus, if $H$ is not \mpcn, then $H$ is not \rcn.

Moreover, from any graph $H$ that is not \mpc it is easy to come up with many graphs
that are not weakly \rcn: Let $H'$ be obtained by attaching disjoint copies of
large complete graphs of slightly different sizes to the vertices $v\in V(H)$ and let $G$ be obtained by doing the same at each corresponding root $v_1$ in $M'(H)$.  The complete graphs ensure that
any $H'$-minor in $G$ is a rooted $H'$-minor; then removing the complete graphs yields a rooted
$H$-minor in $M'(H)$, which is a contradiction since $H$ is not \mpcn.

\medskip

To study the nature of graphs that are not \mpc we now introduce the notion of a shift automorphism. Recall that a bijection $\pi:V\to V$
is called an {\em automorphism} of a  graph $G=(V,E)$ if $uv\in E$ implies
$\pi(u)\pi(v)\in E$.  Recall that for disjoint subsets $A,B$
of vertices in $V$, $[A,B]=\{uv\in E: u\in A,v\in B\}$.

\begin{definition} We say that an automorphism $\pi$ of a
graph $G$ is a {\em shift automorphism} if every vertex
$v$ is adjacent to $\pi(v)$.
\end{definition}

\begin{definition} Let $H$ be a graph and let $S$ be a stable set in $H$.
We say that $S$ {\em induces an \mpcontn} if
 there is a matching in $[S,N(S)]$ that covers all vertices in $N(S)$, and
 $H\graphminus (S\cup N(S))$ has a shift automorphism.
\end{definition}

\begin{prop}\label{M'rout}
If $H$ has a stable set inducing an \mpcontn,
then $H$ is \mpcn.
\end{prop}

\begin{proof} We must show that $M'(H)$ has $H$ as a rooted minor.
Let $\pi$ be the shift automorphism of $G=H\graphminus (S\cup N(S))$.
If $v\in S$, then let $C_v=\{v_1\}$. If $v\in N(S)$, then let
$C_v=\{v_1,v_2,u_2\}$ where $uv$ is the edge covering $v$ in
the matching. Finally if $v\in V(G)$, then let $C_v=\{v_1,\pi(v)_2\}$.

The sets $C_v$ are vertex-disjoint, connected vertex subsets
in $M'(H)$.  It remains to check that if $vw$ is an edge
in $H$, then there is an edge between $C_v$ and $C_w$ in $M'(H)$.
If $v\in N(S)$, then $v_2w_1$ is such an edge. If $vw\in E(H)$
and $v,w\not\in N(S)$, then $v,w\in V(G)$, so $\pi(v)_2\pi(w)_2$
is an edge from $C_v$ to $C_w$.
\end{proof}

We immediately obtain:

\begin{cor}\label{bipM'}
If $H$ is bipartite, then $H$ is \mpcn.
\end{cor}

\begin{proof} Let $W$ be a minimum vertex cover of $H$
and let $S=V(H)\setminus W$.  Then $S$ is a stable set, $N(S)=W$,
and there is
a matching in $[S,W]$ which covers $W$.  $H\graphminus (S\cup W)$ is the
empty graph so it has a trivial shift automorphism.
\end{proof}

Characterizing \mpc graphs seems difficult, since
the property is not closed under taking subgraphs:
For any graph $H$, let $H'$ be obtained
by attaching a leaf to every vertex of $H$, and let $S$ be the set of these leaves.
By Proposition~\ref{M'rout}, $H'$ is \mpcn, but its induced subgraph
$H$ need not be \mpcn.
(Moreover, following Corollary~\ref{K7}, we give an example of a spanning subgraph of
an \mpc graph that is not \mpcn.)

We are more interested in examples of
graphs that are not \mpcn, and thus not \rcn. To this end
we establish the following partial converse of
Proposition~\ref{M'rout}:

\begin{theorem} \label{tri} Let $H$ be a triangle-free graph. Then $H$ is \mpc
if and only if $H$ has a stable set inducing an \mpcontn.
\end{theorem}

\begin{proof} It suffices to show that if $H$ is a rooted minor of $M'(H)$,
then $H$ has a stable set which induces an \mpcontn.
Let $C_v$ be the set in $M'(H)$ which is contracted to
obtain $v$ in the rooted minor; then $v_1\in C_v$. Let
$V^1=\{v: C_v=\{v_1\}\}, V^2=\{v: |C_v|=2\}$,
$V^3=\{v: |C_v|>2\}$, and let $W=V(H)\setminus (V^1\cup N(V^1))$.
Then $\{V^1,V^2,V^3\}$ and $\{V^1,N(V^1),W\}$ are partitions of $V(H)$.
For any $S\subseteq V(H)$, let $S_1=\{v_1: v\in S\}$ and $S_2=\{v_2: v\in S\}$.

For each edge $uv$ in $H$, any $u_1,v_1$-path in $M'(H)$ has length at least three,
because $u_1v_1$ is not an edge, and a path of the form $u_1w_2v_1$ would imply that $uwv$
forms a triangle in $H$. Thus, if $u\in V^1$ and $uv\in E(H)$,
then $v\in V^3$. It follows that $V^1$ is a stable set and $N(V^1)\subseteq V^3$.

Each $v\in N(V^1)$ has a neighbor $u\in V^1$, and $C_v$ must contain a vertex 
adjacent to $C_u=\{u_1\}$, which must be in $N(V^1)_2$.
Then each $C_v$ with $v\in N(V^1)$ contains a distinct vertex of $N(V^1)_2$,
so $\bigcup_{v\in N(V^1)}C_v$ contains $N(V^1)_2$.
Since $\bigcup_{v\in N(V^1)}C_v$ also contains $N(V^1)_1$
and $N(V^1_1\cup V^1_2)\subseteq N(V^1)_1\cup N(V^1)_2$,
there are no edges from $(V^1_1\cup V^1_2)\setminus \bigcup_{v\in N(V^1)}C_v$
to $(W_1\cup W_2) \setminus \bigcup_{v\in N(V^1)}C_v$.
Hence, $C_w\subseteq W_1\cup W_2$ for every $w\in W$.

Since $W\subseteq V^2\cup V^3$, we get
$2|W|\leq \sum_{w\in W}|C_w| \leq |W_1\cup W_2|= 2|W|$, so $W\subseteq V^2$
and $\bigcup_{w\in W}C_w=W_1\cup W_2$.
Since $N(V^1)\subseteq V^3$, we get $W=V^2$ and $N(V^1)=V^3$.

We define a shift automorphism $\pi$ on $H[W]$ (the subgraph of $H$
induced by $W$) as follows: If $w\in W$,
then $C_w=\{w_1,v_2\}$ for some $v\in W$; let $\pi(w)=v$.
Note that $\pi$ is a permutation of $W$.
For each $w\in W$, $w_1\pi(w)_2$ is an edge in $M'(H)$, so
$w\pi(w)$ is an edge in $E(H)$.  If $uv$ is an edge in $H[W]$, then
$M'(H)$ has a $u_1,v_1$-path on $C_u\cup C_v$ of length three or more,
which can only be $u_1\pi(u)_2\pi(v)_2v_1$; hence, $\pi(u)\pi(v)$ is also
an edge.  Thus, $\pi$ is a shift automorphism on $H[W]$.

It remains to show that there is a matching between $V^1$ and $N(V^1)$
that covers $N(V^1)$. Consider any $v\in N(V^1)$. Then there is a vertex
$u\in V_1\cap N(v)$, and
there is a $u_1,v_1$-path $P$ of length at least three on $C_u\cup C_v$.
Since $C_u=\{u_1\}$, the neighbor of $u$ on $P$ is a vertex $v_2'$ in
$C_v\cap N(V^1)_2$.  Since $|N(V^1)_2|=|N(V^1)|$, each vertex $v\in N(V^1)$
has $|C_v\cap N(V^1)_2|=1$, or $C_v\cap N(V^1)_2=\{v_2'\}$.
Also, $C_v\cap N(V^1)_1=\{v_1\}$, $\bigcup_{u\in V^1}C_u=V^1$, and
$\bigcup_{w\in W}C_w=W_1\cup W_2$, so
$C_v\setminus\{v_1,v_2'\}\subseteq V^1_2$.  Therefore, $V^1_2$ contains the
vertex $u_2'$ adjacent to $v_1$ on $P$.  Then $u'v$ is an edge in $H$,
and since $u_2'\in C_v$ and sets of the form $C_v$ are disjoint,
the vertices $u'$ are distinct for different $v$.  Hence,
the edges of the form $u'v$ for all $v\in N(V^1)$, form a matching.
\end{proof}

Suppose that $H$ is a graph with a non-empty stable set $S$ that induces an
\mpcontn, and let $M$ be a matching on $[S,N(S)]$ that covers $N(S)$.
If every vertex in $S$ has degree at least 2,
then $[S,N(S)]$ contains an $M$-alternating cycle:
starting at $v_0\in S$ form a
walk $v_0,v_1,\dots$ such that $v_{2i-1}v_{2i}$ is in the matching
(possible since $N(S)$ is covered) and $v_{2i}v_{2i+1}$ is not
in the matching (possible since every vertex in $S$ is incident to
at least 2 edges). This walk will eventually repeat at a vertex in $S$
and yield the desired cycle.  This observation enables us to prove
the following result.

\begin{cor}\label{2odds}
If $H$ is connected and contains 2 odd cycles of length at least 5 that share
at most one vertex, then $H$ is not \rcn.
\end{cor}

\begin{proof}
It suffices to show that if $H$ is a graph consisting of two such odd cycles
joined by a (possibly  trivial) path, then $H$ is not \mpcn, or equivalently,
$H$ has no stable set $S$ that induces an \mpcontn.
Indeed, such a set $S$ would need to be non-empty, since $H$ has no
shift-automorphism itself, but then the argument above shows that $H$ would
contain an even cycle, a contradiction.
\end{proof}

It immediately follows that $K_9$ and any cactus containing
two long odd cycles are not \rcn.
On the other hand, by combining the ideas of Proposition~\ref{M'rout} and
Lemma~\ref{triangle}, we can show that every cactus not covered by this result is
\mpcn.

\begin{theorem} If a cactus has at most one odd cycle of length at least 5, then it is \mpcn.
\end{theorem}

\begin{proof}
Let $H$ be a cactus with at most one odd cycle of length at least 5.
We will use induction on $|V(H)|$ to construct vertex sets $C_v$
that form a rooted $H$-minor in $M'(H)$.

If $H$ has only one block, then the result follows from Corollary~\ref{cactus}.
Thus we may assume that $H$ has at least two leaf-blocks, one of which is not an odd cycle of length at least 5. This leaf block can be a leaf with its neighbor and incident edge, an even cycle, or a triangle.

If $H$ has a leaf $y$ with neighbor $x$, apply induction to find a rooted minor of $H\graphminus\{x,y\}$ in $M'(H\graphminus\{x,y\})$,
then let $C_y=\{y_1\}$ and $C_x=\{x_1,x_2,y_2\}$.  Since $x_2$ is adjacent to $z_1$ for every $z\in N_H(x)$, and
$y_2$ is adjacent to $x_1$ and $x_2$, this suffices.

If $H$ has a triangle leaf-block with vertices $x,y,z$ such that $z$ is the cut-vertex, apply induction
to find a rooted minor of $H\graphminus\{x,y\}$ in $M'(H\graphminus\{x,y\})$, then let $C_x=\{x_1,y_2\}$ and $C_y=\{y_1,x_2\}$.  

If $H$ has an even cycle leaf-block $B$ with vertices labeled $x^0,x^1,\ldots,x^{2k}=x^0$ with cut-vertex $x^0$,
apply induction to find a rooted minor of $H\graphminus V(B)$ in $M'(H\graphminus V(B))$,
let $C_{2i}=\{x^{2i}_1,x^{2i}_2,x^{2i+1}_2\}$ and $C_{2i+1}=\{x^{2i+1}_1\}$ for $0\le i<k$.
\end{proof}

The cycles in a  cactus have a very simple structure. From this point of view, the next graphs to consider
are the Theta graphs. For $k,l,m\ge 0$, at most one of which is zero, the {\em Theta graph}
$\theta(k,l,m)$ is the $(2+k+l+m)$-vertex graph obtained
by connecting two vertices $x,y$ by three paths which contain
$k,l,m$ internal vertices, respectively. For example,
$K_{2,3}=\theta(1,1,1)$.

\begin{theorem}\label{theta}
The Theta graph $H=\theta(k,l,m)$ is not \mpc
if and only if exactly one of $k,l,m$ is odd and $H$ is triangle-free. 
\end{theorem}

\begin{proof}
If $k,l,m$ all have the same
parity, then $H$ is bipartite, so by Corollary~\ref{bipM'}, $H$ is \mpcn.
Suppose that exactly one of $k,l,m$ is even. Let $z$ be
the neighbor of $x$ on the unique $x,y$-path of odd length. Then $H\graphminus xz$ has a unique
bipartition $[S,N(S)]$ with $x,y,z\in N(S)$, with a matching that covers $N(S)$.
By Proposition~\ref{M'rout}, $H$ is \mpcn.
Now suppose that exactly one of $k,\ell,m$ is odd and $H=\theta(k,l,m)$ contains a triangle.
Then $H$ consists of an odd cycle $v^1,v^2,\dots,v^n,v^1$ plus
the edge $v^1v^3$. Letting $C_{v^i}=\{v^i_1,v^{i+1}_2\}$ yields a rooted $H$-minor.
This finishes the proof of the forward direction.

We prove the backward direction by contradiction, so
suppose that exactly one of $k,l,m$ is odd and $H$ is triangle-free,
and let $S$ be a stable set in $H$ that induces an \mpcontn.
If $S=\emptyset$, then $H$ must have a shift isomorphism $\pi$.
Since $x$ and $y$ are the only vertices of degree 3,
 $\pi(x)=y$ and $\pi(y)=x$, and $xy\in E(H)$.
 Since $H\graphminus \{x,y\}$ is acyclic,
the permutation $\pi$ has no cycles with three or more elements.
Then $\{v\pi(v): v\in V(H)\}$ is a perfect matching in $H$,
so $|V(H)|$ is even, a contradiction.

Thus, we may assume that $S$ is non-empty.  Since $S$ induces
an \mpcontn, we can let $M$ be the matching in
$[S,N(S)]$ that covers $N(S)$. Because every vertex has degree
at least 2, it follows that $[S,N(S)]$ contains an $M$-alternating
cycle $C$. Since $C$ is an even cycle it must consist of the two $x,y$-paths of odd length.
Since their lengths are odd, only one of $x,y$
is in $N(S)$, say, $x\in S$. Then $x$ has a unique neighbor $z$ in $N(S)\setminus V(C)$.
Since $M$ covers $N(S)$, we can construct an $M$-alternating path beginning with $x,z,\dots$,
which shows that every edge in $H$ is in $[S,N(S)]$.
Therefore $H$ is bipartite, a contradiction.
\end{proof}

This yields two graphs on 7 vertices which are not \mpcn,
and therefore not \rcn: $\theta(0,2,3)$ (a 7-cycle with
a chord) and $\theta(1,2,2)$ (a 6-cycle with a subdivided diagonal).
Also any graph containing $\theta(0,2,3)$ or $\theta(1,2,2)$ is not \rcn,
by Lemma~\ref{subgraph}. 
It can be checked by a tedious case analysis that every graph $H$ on at most 6 vertices
is \mpcn.

We now combine the results from this section to prove a simple necessary condition for a graph to be \rcn.
A graph is {\em $k$-critical} if it has chromatic number $k$, but every one of its proper subgraphs
is $(k-1)$-colorable. A $k$-critical graph for $k>2$ must be 2-connected.
We will call an odd cycle on at least 5 vertices {\em long}.

\begin{theorem}\label{chromatic bound}
If $G$ is contractible, then $G$ is 6-colorable.
\end{theorem}

\begin{proof}
We proceed by contrapositive, so suppose $G$ is not 6-colorable.
By Lemma~\ref{subgraph} we may assume that $G$ is 7-critical, and we can consider a specific
7-coloring of $G$.

We first show that any 4 color classes in $G$ must contain a subgraph that is a long odd cycle or a $K_4$. Let $H$ be a 4-critical
subgraph of these 4 color classes, and suppose that $H$ contains no long odd cycle. Then $H$ contains a triangle $T$, and for every vertex $v$ in $V(H)-T$ there are two paths from $v$ to $T$. In order to avoid a long odd cycle these paths must be single edges,
so $v$ has at least 2 neighbors on $T$. If two vertices in $V(H)-T$ have different neighbors in $T$ then we obtain a 5-cycle, so
$H$ contains a spanning $K_{1,1,t}$. Moreover, since $H$ is not 3-colorable there is an edge not incident with a vertex in $T$,
and it follows that $H=K_4$.

So let $H$ be the subgraph of some 4 color classes that is either a long odd cycle $C$ or a $K_4$.
The remaining 3 color classes must contain an additional odd cycle $C'$.
Since $G$ is 2-connected there are two vertex-disjoint paths connecting $C'$ and $H$, and we let $x,y$ be the vertices of
these paths that are in $H$. In either case we see that $G$ contains a triangle-free Theta graph that is not contractible by Theorem~\ref{theta}:
 if $H=K_4$ we add an odd length $x,y$-path via $C'$ to obtain $\theta(0,2,2m+1)$, and if $H=C$, we add
an even length $x,y$-path via $C'$ to $H$. Thus the result follows by Lemma~\ref{subgraph}.
\end{proof}

With more effort this result can certainly be improved. For example, it can be shown that $K_6$ is the only 6-critical graph that may be \rcn.

\begin{cor}\label{K7}
 $K_7$ is \mpcn, but not \rcn.
\end{cor}

\begin{proof} $K_7$ has many shift-automorphisms, but is not 6-colorable.
\end{proof}

The result that $K_7$ is not \rc was also obtained independently by
Melody Chan and Paul Seymour~\cite{Se}.

Theorem~\ref{theta} and Corollary~\ref{K7} show that a subgraph of an \mpc graph need not be \mpcn, so there is probably no simple characterization of \mpc graphs.
On the other hand, every proper connected subgraph of a Theta graph is a cactus, so
by Corollary~\ref{cactus} every non-\rc Theta graph mentioned in Theorem~\ref{theta}
is part of the forbidden subgraph characterization of \rc graphs.

\section{Open questions}\label{question-sec}

\begin{enumerate}
\item Is a subgraph of a weakly \rc graph weakly \rcn?
\item Are $K_{2,4}$ or $K_{3,3}$ \rcn?
\item  Are bipartite Theta graphs \rcn?
\item Find a bipartite graph that is not (weakly) \rcn.
\item Are $K_5$ or $K_6$ \rcn?
\item Is every $K_n$ weakly \rcn? How is this question related to Hadwiger's conjecture?
\item Find the forbidden subgraph characterization for \rc graphs.
\item ({\em Decidability}) Given a graph $H$ is there a finite process for deciding if $H$ is (weakly) \rcn?  This is not a given,
since the paths in an \Shsch can get arbitrarily long, and since there are infinitely many forbidden subgraphs in the characterization.
\end{enumerate}

\section*{Acknowledgements}
We thank Paul Seymour and Paul Wollan for useful discussions about
this topic.


\end{document}